\input amstex
\documentstyle{amsppt}
\NoBlackBoxes

\def\st{:}

\def\E{{\Cal E}}

\def\eg{{\it e.g.,}}

\def\rn{\Bbb R^{n}}
\def\rr{\Bbb R}
\def\rp{{\Bbb R_+}}

\def\sxor{{S(x_0,r_0)}}
\def\txo{T_{x_0}}

\define \Dl{\Delta}
\define \dl{\delta}

\define \s{\sigma}
\define \Si{\Sigma}
\define \th{\theta}

\define \vp{\varphi}

\def\capl{\operatornamewithlimits{\bigcap}\limits}

\define \fc{\frac}
\define \iy{\infty}
\define \la{\langle}
\define \ra{\rangle}
\define\pf{\demo{Proof}}
\define \edm{\qed\enddemo}
\define \ep{\endproclaim}

\define \bk{\bigskip}

\define \BR{\Bbb R}

\define \Cd{\Cal D}
\define \CE{\Cal E}

\define \fin{\operatorname{fin}}
\define \supp{\operatorname{supp}}
\define \codim{\operatorname{codim}}

\define \rad{\operatorname{rad}}

\define \spa{\operatorname{span\ }}

\leftheadtext{M.L.\ Agranovsky and E.T.\ Quinto}
\rightheadtext{Geometry of Stationary Sets for the Wave Equation}

\topmatter

\title{Geometry of Stationary Sets for the Wave Equation in $\rn$.\\
 The Case of Finitely Supported Initial Data\\An Announcement}\endtitle

\author Mark L. Agranovsky$^*$ and Eric Todd Quinto$^{**}$\endauthor

\affil Bar Ilan University and Tufts University\endaffil

\address Department of Mathematics, Bar Ilan University, 53000 Ramat
Gan, Israel\endaddress
\email agranovs\@macs.cs.biu.ac.il \endemail

\address Department of Mathematics, Tufts University, Medford, MA 02155
USA \endaddress
\email equinto\@math.tufts.edu\endemail

 \abstract We consider the Cauchy problem for the wave equation in
the whole space $\BR^n,$ with initial data which are distributions
supported on finite sets.  The main result is a precise description
of the geometry of the sets of stationary points of the solutions to
the wave equation.\endabstract

\thanks $^*$ Partially supported by  Grant 408/97-3 from the Academy of
Sciences of Israel.  \endgraf $^{**}$ Partially supported by US
National Science Foundation grants 9622947 and 9877155. \endthanks



\endtopmatter

\document

\subheading{\S1.\ Introduction}

Our goal is to understand the structure of stationary (nodal) sets of
solutions to the wave equation $\BR^n,\quad n \ge 2:$ $$\aligned
&u_{tt}=\Dl u,\ u=u(x,t),\ x\in\BR^n,\ t>0,\\ &u\big|_{t=+0}=0,\quad
u_t\big|_{t=+0}=f\, .\endaligned\tag1.1$$ We will consider the case $f\in
\E'_{\fin}(\rn)$, the set of distributions supported on a finite set of
points.  Stationary sets are sets of points $x\in \rn$ for which the
solution to the wave equation is always zero.  This article is an
announcement of the results in \cite{AQ4}; complete proofs are given there.

We use standard notation for function spaces; for example,
$\Cd(\BR^n)$ is the space of all $C^\iy$-functions with compact
support; $\CE'(\BR^n)$ is the space of compactly supported
distributions; $C_{\rad}^\iy(\BR^n)$ and $\Cd_{\rad}(\BR^n)$ are the
subspaces of corresponding spaces, consisting of radial functions $f,
i.e., f(x)=f(|x|)$.  Finally, let $R_+=(0,\infty)$.

For classical solutions to (1.1), we define stationary sets as follows.

\definition{Definition 1} Let $f\in C^\iy(\BR^n)$ and $u$ be the
(classical) solution for (1.1).  Define the stationary set $S(f)$ as
the set of time-invariant zeros of the solution $u:$ $$S(f)=\{x\in
\BR^n\st u(x,t)=0,\ t>0\}.\tag1.2$$
\enddefinition

 We use regularization to extend Definition 1 to distributional
solu\-tions.  Namely, if $f\in \Cd'(\BR^n)$ and $\vp\in
\Cd_{\rad}(\BR^n)$, the convolution $f\ast\vp$ is smooth and $u\ast
\vp$ (convolution with respect to $x)$ is in $C^\iy$ and solves (1.1)
for the data $f\ast\vp.$

\definition{Definition 2}
For $f\in \Cd'(\BR^n)$ define $$S(f)=\capl_{\vp\in
\Cd_{\rad}(\BR^n)}S(f\ast\vp),\tag1.3$$ where $S(f\ast\vp)$ is defined
by (1.2).\enddefinition

For continuous $f$, one can use Theorem 2 below to show Definition 2
coincides with Definition 1.  In fact, the set of all $\vp\in
\Cd_{\rad}(\BR^n)$ in (1.3) can be replaced by an $\dl$-sequence
$\vp_n,$\ $\vp_n\in
\Cd_{\rad}(\BR^n).$

The main question under consideration is the following.

\proclaim{Problem} Which sets $S\subset \BR^n$ are stationary sets,
$S=S(f),$ for some $f\in \Cd'(\BR^n)$\ $(f\in\CE'(\BR^n))?$\ep

The problem has been solved in the plane \cite{AQ1}, \cite{AQ2} if
$f$ is an arbitrary distribution of compact support, but not much is
known in general.  For distributions in the plane, the stationary
sets have very restrictive structure; they must be the union of a
finite set and a Coxeter system of lines (lines through one point
generated by a finite rotation group).  This Coxeter set is contained
in a translate of the zero set of a homogeneous harmonic polynomial,
and it is conical about the point of intersection.  Loosely speaking,
we will prove that a similar pattern occurs for finitely supported
distributions in $\rn$.

Characterizing stationary sets in $\rn$ for $n>2$ is more difficult
and only partial results are known.  It is known for compactly
supported initial data in $\rn$ that stationary sets are contained in
zero sets of harmonic polynomials union algebraic varieties of lower
dimension, and it is conjectured that the harmonic polynomial can be
assumed to be a translate of a homogeneous polynomial \cite{AQ1},
\cite{AQ2}.  It is shown in \cite{ABK} for $f$ sufficiently
integrable at infinity that stationary sets cannot have bounded
closed components.  In \cite{A} and \cite{AVZ} more precise analyses
are given in $\rn$ for stationary sets of lower dimension and conical
stationary sets for $f$ with arbitrary growth.

In this article we describe stationary sets for the case of the
initial data with finite support for arbitrary dimension $n.$ We
prove that, up to a low-dim\-en\-sion\-al component, the stationary
sets are affine algebraic cones with a special geometry.

A {\it cone} is understood to be a union of straight lines with a
common point which is the vertex of the cone.  We will call a cone
$K\subset \BR^n$ \ {\it $k-$flat} with edge $L,$ where $L$ is a
$k$-dimensional plane in $\BR^n,$ if $K$ is a union of $(k+1)$-planes
containing $L.$

A union $\Sigma=H_1\cup\dots\cup H_q$ of hyperplanes
$H_i\subset\BR^n$ is called a {\it Coxeter system} of hyperplanes if
$\Sigma$ is invariant with respect to any reflection $\s_i$ around
the hyperplane $H_i,$\ $i=1,\dots,q.$ The {\it Coxeter group}
generated by the reflection $\s_1,\dots,\s_q$ will be denoted by
$W(\Si).$

We will call a polynomial $P$ in $\BR^n,$ with real coefficients {\it
a harmonic divisor} if $P$ divides a nonzero harmonic polynomial.
Zero sets of homogeneous harmonic divisors will be called {\it
harmonic cones}.

For any set $F\subset \BR^n,$ the affine subspace spanned by $F$ will
be denoted by $\spa F.$

\medskip
Our main result is the following.  \proclaim{Theorem 1} Let
$f\in\CE_{\fin}'(\BR^n)$,\ $f\ne0.$ If $S(f)\ne\emptyset,$ then
\roster\item"(a)" $S(f)$ is an algebraic variety in $\BR^n$,
contained in the zero set of a nonzero harmonic polynomial.
\item"(b)" After a suitable translation, the set $S(f)$ can be
represented in the form $$S(f)=S_0\cup V,$$ where $V$ is an algebraic
variety of $\codim V>1$ and $S_0,$ assuming it is nonempty, is a
harmonic cone, which is a $(n-1)-$dimensional real algebraic variety.

\endroster 
In addition, the following is true:
\roster\item"(c)" The conical component $S_0$, in general, has the
two components $$S_0=\Si\cup K,$$ where $K$ contains $\supp f$ but
$\Si$ does not, $\Si\cap K\ne\emptyset$ provided both $\Si$ and $K$
are nonempty, $\Si$ is a Coxeter system and $K$ is a $k$-flat
harmonic cone with  the edge $L=\spa(\supp f),$\ $k=\dim L<n.$ If
$\supp f$ is a generic set, i.e., $k=n,$ then $K=\emptyset.$ If
$k=n-1$, then $K$ is a hyperplane and $\Si\cup K$ is a Coxeter
system.  
\item"(d)" If $\tilde\Si$ is the union of all hyperplanes
contained in $S_0,$ then $\tilde \Si$ is again
a Coxeter system; the
distribution $f$ is odd with respect to any reflection $\s\in
W(\tilde\Si),$ i.e. $f\circ \s=-f;$ the sets $S_0,$\ $V,$\ $S(f)$
and $\supp f$ are $W(\tilde\Si)$-invariant.
\endroster\ep

Theorem 1 says that finite sets of point sources generate stationary
sets which are necessarily algebraic varieties and which are either
small (empty or low-dim\-en\-sion\-al) or up to a low-dimensional
component, are $(n-1)$-dimensional cones which suitably translated
are determined by zeros of spatial harmonics.

The geometry of the essential, conical part is as follows.  If the
set of points in $\supp f$ is generic then the cone is a Coxeter
system of hyperplanes.  These stationary sets may appear only as a
result of a Coxeter skew-symmetry of the initial data.  If $\supp f$
lies in a proper affine subspace in $\BR^n$, then another component
may appear which is a cone containing $\supp f$.  In the plane, for
any compactly supported $f$, $S(f)$ is, up to a finite set, 
a Coxeter system of lines \cite{AQ1}.
However, in the plane the set $K$ in Theorem 1 would be a collection
of lines and therefore a Coxeter system by (d).

An important problem in studying of the wave equation is
characterizing {\it nodal} sets (see \cite{CH}, I, Ch.5, S.5), that is
zero sets of eigenfunctions of the Laplace operator, or,
equivalently, zero sets of time-harmonic solutions of the wave
equation. This problem has been studied by many authors.  Results on
this subject mainly say that nodal sets are hypermanifolds with
singularities and the eigenfunctions cannot vanish to high order on
the nodal sets (see,\eg\ \cite{DF1}, \cite{DF2}, \cite{Ch}, \cite{B1}
and others).

The problem under consideration is directly related to describing
nodal sets. Indeed, extending the solution $u(x,t)$ of (1.1) for $t <
0$ by $u(x,t)=-u(x,-t)$ and applying Fourier transform in $t$ to the
both sides of (1.1) yields $$ -\lambda^2 v(x, \lambda)= \Dl
v(x,\lambda),$$ where $v(x,\lambda)$ is the Fourier transform
evaluated at arbitrary $\lambda \in\BR$.

Thus the stationary set $S(f)$ (1.2) is just the intersection of
nodal sets of all the eigenfunctions $v(\cdot ,\lambda)$ which are,
since the initial data $f$ has compact support, nonzero for an
infinite number of $\lambda$.  Thus, while for a single eigenfunction
or a finite linear combination of eigenfunctions, the available
information is the general analytic structure of nodal sets, our
result shows that joint nodal hypermanifolds of one-parameter
families of eigenfunctions are strongly determined geometrically.  It
is worth noting in this connection the result of \cite{B2} which
describes geometry of nodal lines of bounded membrane under the
assumption that these lines contain a open piece of a straight line.

This announcement continues a series of works \cite{AQ1}, \cite{AQ2},
\cite{AQ3}, \cite{ABK}, \cite{A}, \cite{AVZ}, \cite{AR}, started by
\cite{AQ1} and devoted to the description of injectivity sets for the
spherical transform, stationary sets for the wave and heat equations,
and related problems.  Our initial interest in the problem was
motivated by a problem in approximation theory posed in \cite{LP}
(cf. \cite{AQ1}, \cite{AQ2}).  Complete proofs of our results are
given in \cite{AQ4}.

\bigskip
\subheading{\S2.\ Proof outline for Theorem 1}

Let $f\in C(\BR^n)$.  Define the spherical transform $C(\BR^n)\ni
f\to\hat f\in C(\BR^{n}\times\rp)$ by $$\hat f(x,r)=Rf(x,r)=
\int_{|\th|=1}f(x+ r\th)dA(\th),\tag2.1$$ where $dA$ is the
normalized surface measure on the unit sphere.  Because this
spherical transform is a Fourier integral operator, it can be defined
on distributions.

 \proclaim{Theorem 2} Let $f\in C(\BR^n)$ and let $S(f)$ be defined by
$(1.2).$ Then $$S(f)=\{x\in \BR^n\st \hat f(x,r)=0\quad\text{for all}\quad
r>0\}\tag 2.2$$ and also $$S(f)=\{x\in\BR^n\st (f\ast u)(x)=0\quad\text{for
all}\quad u\in\Cd_{\rad}(\BR^n)\}.$$\ep

The second equality is a key to proving the consistency of
Definitions 1 and 2 for smooth functions.  The proof of Theorem 2
follows from the Poisson-Kirchoff formula and integral equations
techniques (e.g., \cite{Q}).

For $f\in \CE'(\rn)$, one can show that $S(f)$ satisfies equality
(2.2) distributionally.  One
needs to observe that, for each $x\in \rn$ one can define $\hat f(x,\cdot)$
as a distribution on even functions in $\Cd(\rr)$ using the natural map
between this set of functions and the set of smooth, compactly 
supported radial functions centered at $x$.

If $Q$ is a polynomial, we let $N(Q)$ denote the zero set of $Q$ in
$\rn$.  Using ideas related to those in \cite{AQ1}, \cite{AQ2}, we
prove

\proclaim{Theorem 3} Let $f\in\CE'(\BR^n).$ Then
\roster\item"a)" The set $S(f)$ is an algebraic variety in $\BR^n$
contained in the zero sets of a nonzero harmonic polynomial.
\item"b)" $S(f)=S_0 \cup V,$ where $V$ is an algebraic variety of
$\codim V>1$, $S_0=N(Q)$  and $Q$ is a harmonic divisor.\endroster\ep

This gives us the key to understand the algebraic structure of $S(f)$.

We now need to understand the geometric structure of $S(f)$.  We
introduce some new notation.  Let $r\in \rp$, $S\subset \rn$, and
$x\in S$.  The point $x$ is called a {\it regular point} of $S$ if
and only if there is a connected real-analytic hypersurface, $A$, (an
$(n-1)-$dimensional real-analytic submanifold of $\rn$) such that
$x\in A\subset S$.  Let $x$ be a regular point of $S$, and let $A$ be
such an associated hypersurface ($x\in A\subset S$).  Then, we let
$T_x$ denote the hyperplane tangent to $A$ at $x$.  The points $y$
and $y'$ in $\rn$ are said to be {\it $T_x-$mirror} if and only if
they are reflections about $T_x$.  If $y \in T_x$, then $y$ is its
own mirror point, and we say $y$ is {\it self-mirror}.  Our next
theorem is a microlocal version of a reflection principle of Courant
and Hilbert.

\proclaim{Theorem 4} Let $f\in\E'_{\fin}(\BR^n)$,\ $f\ne0$.
We let $x_0 \in S=S(f)$ be a regular point.  Let $ A\subset S$ be a
connected real-analytic hypersurface containing $x_0$, and let $\txo$
be the hyperplane tangent to $A$ at $x_0$.  Let $y_0\in \supp f
\setminus {x_0}$.  Then, the $\txo-$mirror point to $y_0$ must also
be in $\supp f$.\endproclaim

Courant and Hilbert \cite{CH {\rm II, pp. 699 {\it ff.}}} proved the
{\it reflection principle} that if $A=\txo$ is a hyperplane and $f$
had zero integrals over all spheres centered on $\txo$, then $f$
would have to be an odd function about $\txo$.  Therefore, if $f$ had
zero integrals over such spheres, and $f$ were zero at one mirror
point, it would have to be zero at the other.  This theorem is a microlocal
version of that fact.  It is proven by calculating the microlocal
properties of the spherical transform and using a theorem of Kawai,
Kashiwara and H\"ormander about analytic wavefront at boundary points
of supp~$f$ \cite{H\"o}.

 \demo{Proof sketch for Theorem 1} Let $f\in \E'_{\fin}(\rn)$.  Let
$x_0\in A\subset S$ be a regular point, let $A$ be a real-analytic
hypersurface such that $x_0\in A\subset S$.  Let $\txo$ be the plane
tangent to $A$ at $x_0$.  Using the result of Theorem 3, we have
$S(f) = N(Q)\cup V$ where $Q$ is a harmonic divisor; $S_0=N(Q)$.

The geometric part of the proof has two cases.

First, assume $\supp f\not\subset \txo$.  Then, there is an $r_0>0$
and a point $a_0 \in \sxor$ that is in $\supp f\setminus \txo$.  By
the support theorem, the mirror point, $a_m$, to $a_0$ must also be
in $\supp f$.  By definition, $\txo$ is the perpendicular bisector of
the segment $\overline{a_0a_m}$.  For $x\in A$ near $x_0$, the sphere
$S(x,|x-a_0|))$ must meet $\supp f$ at $a_0$ and at its $T_x$ mirror
point by Theorem 4.  Because there are only a finite number of
points in $\supp f$ and $x$ is close to $x_0$, this mirror point must
be $a_m$.  So, $T_x$ is the perpendicular bisector of
$\overline{a_0a_m}$ for all $x\in A$ sufficiently close to $x_0$.
Therefore, $T_x=\txo$.  This shows locally near $x_0$, $A$ is flat.
Because $S(f)$ is an algebraic variety, this shows that $\txo\subset
S(f)$.  Now, we use the reflection principle \cite{CH {\rm II, pp.
699 {\it ff.}}} to conclude $f$ is odd about $\txo$.  

We do this construction for each regular point in $x\in S$ for which $\supp
f \not\in T_x$.  This gives us a collection of hyperplanes contained in
$S(f)$ such that $f$ is odd about each one.  We know this collection is
finite since $S(f)$ is an algebraic variety.  It also gives us a Coxeter
group of reflections about these hyperplanes, and the set of all such
hyperplanes generated by the Coxeter group is $\Sigma$.  Using properties
of Coxeter groups 
\cite{DS, Ch. 8, S.10, Th. 8}, \cite{GB, Prop. 4.1.3}, 
we
show the Coxeter group is finite and the intersection $\cap\Sigma$ is
nonempty.  We use Theorem 3 and the irreducible factors of the harmonic
divisor $Q$ corresponding to $\Sigma$ to show these factors are all
homogeneous about the points in $\cap\Sigma$.

Now, for the second case, assume $\supp f\subset \txo$.  This implies
$\dim L < n$.  Let $a_0\in \supp f$.  In this case, $a_0$ is a
$\txo-$self-mirror point.  If there were an $x\in A$ near $x_0$ for
which $a_0$ was not $T_x-$self-mirror, then this would contradict 
Theorem 4 because by finiteness of $\supp f$, there would be no
$T_x-$mirror point to $a_0$ in $\supp f$.  This shows that locally
$A$ generates a subset of $S_0$ that is conical about $a_0$.  By a
linearity argument, $A$ generates a subset of $S_0$ conical about any
point in $L=\spa \supp f$.  We use Theorem 3 and the irreducible
factors of the harmonic divisor $Q$ corresponding to $K$ to show
these factors are all homogeneous about points in $L=\spa\supp f$.

This finishes the geometric part of the proof.  If $K\neq \emptyset$ and
$\Sigma\neq \emptyset$, we show, using Kakutani theorem and 
invariance of $L$ under $\Sigma-$ reflections, that 
$L$ and $(\cap\Sigma)$ have common points and therefore $S_0=K \cup \Sigma$
is a cone about any such a point.  This finishes the proof
sketch.\edm

\bk

\subheading{\S3.\ Sufficient Conditions for Stationary Sets}

Theorem 1 says that the stationary sets $S(f)$ for
$f\in\CE_{\fin}'(\BR^n)$ may consist of three parts: a
low-dimension\-al variety $V,$ a Coxeter system $\Si$, and a cone $K$
having all points in $\spa(\supp f)$ as vertices.  In addition, the
union $S_0=\Si\cup K$ must be a cone containing the zero set of some
shifted harmonic homogeneous polynomial and the entire stationary set
$S(f)=\Si\cup K\cup V$ must belong to the zero set of some nonzero
harmonic (not necessarily homogeneous, if $V$ is not a cone with the
common vertex with $S_0).$

Now the question is whether all the possibilities are realizable.
Namely, whe\-ther each of the sets $\Si,$\ $K,$\ $V$ and any unions of
the sets of these three types are the stationary set $S(f)$ for some
$f\in\CE_{\fin}'(\BR^n)$ or, more generally, are contained in a
stationary set?

Below we give positive answers for the sets $\Si$, $K$, $V$ and
$\Si\cup V$.  The case of $\Si\cup K\cup V,$ where each of the three
sets are nonempty, remains unsolved.

Given a polynomial $G\in\BR[x_1,\dots,x_n],$ denote by $T_G$ the
distribution $\la T_G,\vp\ra=G(\partial)\vp(0),\ \vp\in \Cd(\BR^n).$ Here
$G(\partial)=G\left(\fc{\partial}{\partial
x_1},\dots,\fc{\partial}{\partial x_n}\right).$

The following theorem shows that the stationary set generated by
a homogeneous distribution (of finite order) supported at a single 
point coincides with common zeros of iterated Laplacians of the
symbol of the corresponding differential operator:

\proclaim{Theorem 5}
For any homogeneous polynomial $G\in\BR[x_1\dots x_n]$, we have
$$S(T_G)=\capl_{j\ge0}N(\Dl^j(G))\, .$$\ep

The following two corollaries prove that zero sets of harmonics are 
stationary sets of a homogeneous distribution supported at a single point
and describe all such distributions:

\proclaim{Corollary 6} If $\Psi$ is a homogeneous harmonic
polynomial, then $N(\Psi)=S(T_\Psi).$\ep

\proclaim{Corollary 7} Let $\Psi$ be a homogeneous harmonic polynomial
and $G$ a polynomial in $\BR^n.$ Then $N(\Psi)\subset S(T_G)$ if and
only if $\Psi$ divides all the polynomials $G,$ $\Delta G,$
$\Delta^2G,\dots\ .$
\ep

\pf Since $\Psi$ is homogeneous, then it is easy to check that
$N(\Psi)\subset S(T_G)$ is equivalent to $N(\Psi)\subset S(T_{G_{m}}),$
where $G_m$ is any homogeneous term of $G.$ In turn, by Theorem 5, this
is equivalent to $G_m,\Dl G_m,\Dl^2G_m,\dots$ vanishing on $N(\Psi).$
It can be proven that
vanishing on zeros of a real harmonic polynomial is equivalent to divisibility
, therefore all the homogeneous terms
$G_m,$ along with their iterated Laplacians, are divisible by $\Psi.$ This
proves the corollary.
\edm

Finally, any low-dimensional real algebraic variety can be stationary
for some solution of the wave equation with point supported initial
data and, moreover, any Coxeter system can be added:
\proclaim{Theorem 8} \cite{A}\ Let $V$ be an algebraic variety in
$\BR^n,$ \ $\codim V>1.$ Let $\Si$ be either empty or a Coxeter
system of hyperplanes.  Then there exists a nontrivial polynomial
$G\in\BR[x_1,\dots,x_n]$ such that $\Si\cup V\subset S(T_G).$\ep

\remark{Remark} In order to prove that $N(\Psi)\cup V,$\ $\Psi$ is a
homogeneous harmonic polynomial, $\codim V>1,$ can be realized as a
stationary set, it would be sufficient to prove, according to
Corollary 7, that the set of homogeneous polynomials $G$ such that $\Psi$
divides all $\Dl^sG,$\ $s=0,1,\dots$ is big enough to satisfy the
additional condition on the low-dimensional part:
$\Dl^sG\big|_V=0,$\ $s=0,1,\dots\ .$ Due to what has been
proven in Theorem 5, this means that all the harmonic homogeneous
polynomials $h_{k-2j}$ in the decomposition 
$$G(x)=h_k(x)+|x|^2h_{k-2}(x)+..., \quad k=\deg G$$
are divisible by
$\Psi$ and vanish on $V.$ However, the question about whether the
space of harmonic homogeneous polynomials $h$ divisible by a given
harmonic $\Psi$ is big enough turned out to be very nontrivial in
$\BR^n$ for $n>2$ ({\it cf.} \cite{A}). We even do not know whether
this space is always infinite dimensional or not.\endremark

\subheading{\S4.\ The Case of Balls}

Similar arguments can be used to prove a theorem similar to Theorem 1
if $\supp f$ is the disjoint union of balls.  Let $\CE_D'(\rn)$ be
the set of distributions whose support is the disjoint union of a
finite number of closed balls.  Here is the geometric analogue of
Theorem 1 for $\CE_D'(\rn)$.

 \proclaim{Theorem 9\ (Support Theorem)}\ Let
$f\in\CE_{D}'(\BR^n)$,\ $f\ne0.$ Assume $S=S(f)\ne\emptyset$.  Assume
there are regular points in $S$, and let $x_0\in S$ be a regular point
and $x_0\notin\supp f$.  Let $A$ be a connected real-analytic
hypersurface in $\rn$ such that $x_0\in A\subset S$.  Let $\txo$ be the
hyperplane tangent to $A$ at $x_0$.  Let $C$ be the set of centers of
the disks making up $\supp f$.  There are two possibilities.

\roster
\item"(a)" For some $c_0\in C$, $c_0\notin \txo$.  In this
case $A\subset \txo\subset S$ and $\supp f$ is symmetric about $\txo$.
Furthermore, $f$ is odd about $\txo$.  

 \item"(b)" Or, $C\subset \txo$.  In this case, near $x_0$, $S$ is conical
about $L=\spa C$.  Precisely, $A$ generates a subset of $S$ that is conical
with edge $L$.  In this case, $k=\dim L <n$.
\endroster\ep

\subheading{\S5.\ Concluding Remarks}

Theorem 1 asserts that for an initial distributions with finite
support, the essential $(n-1)-$dimensional part of the stationary set
is a cone.  From Section 3, we learn this cone appears as the set of
common zeros of spatial harmonics in the Fourier decomposition of the
initial distribution. Correspondingly, this happens only when these
harmonics have a large set of common zeros (are coherent).  More
specifically, the cone may contain a system of Coxeter mirrors, if
the initial data (sources) admit a corresponding symmetry. In this
case vanishing of the solution of the wave equation on the mirrors is
the result of cancelling of waves propagated by symmetric sources.

We expect that the stationary sets have a similar geometry for compactly
supported initial data and, more generally, for distributions vanishing
sufficiently fast at infinity.  The main difficulty in proving that is
obtaining the conical structure of the essential part of stationary sets.
This was done in \cite{AQ1}, \cite{AQ2} for $n=2$. There, the simple
structure of zero sets of harmonic polynomials of two variables and the
support theorem, Theorem 4, play important roles.  Lack of information
about zero sets of harmonic polynomials of more than two variables
was our
main obstacle in extending our approach to $n>2$.  Nevertheless, we hope to
succeed using a deeper analysis of the algebraic and geometric structure of
stationary sets and by refining the microlocal results that go into the
proof of Theorem 1 to be valid more generally, such as for rapidly
decreasing functions.

\bk

\enddocument

\bye
\end